\documentclass{amsart}
\usepackage{amssymb}
\usepackage{amsfonts}
\usepackage{graphicx}
\usepackage{psfrag}
\usepackage{verbatim}
\usepackage{hyperref}

\newcommand{\cellsize}{11}
\newlength{\cellsz} 
\newsavebox{\cell}
\sbox{\cell}{\begin{picture}(\cellsize,\cellsize)
\put(0,0){\line(1,0){\cellsize}}
\put(0,0){\line(0,1){\cellsize}}
\put(\cellsize,0){\line(0,1){\cellsize}}
\put(0,\cellsize){\line(1,0){\cellsize}}
\end{picture}}
\newcommand\cellify[1]{\def\thearg{#1}\def\nothing{}%
\ifx\thearg\nothing
\vrule width0pt height\cellsz depth0pt\else
\hbox to 0pt{\usebox{\cell} \hss}\fi%
\vbox to \cellsz{
\vss
\hbox to \cellsz{\hss$#1$\hss}
\vss}}
\newcommand\tableau[1]{\vcenter{\vbox{\let\\\cr
\baselineskip -16000pt \lineskiplimit 16000pt \lineskip 0pt
\ialign{&\cellify{##}\cr#1\crcr}}}}
\newcommand\tabl[1]{\vtop{\let\\\cr
\baselineskip -16000pt \lineskiplimit 16000pt \lineskip 0pt
\ialign{&\cellify{##}\cr#1\crcr}}}
\newtheorem{lemma}{Lemma}
\newtheorem{thm}{Theorem}
\newtheorem*{thmp}{Theorem $1'$}
\newtheorem*{problem}{Problem}
\newtheorem{prop}{Proposition}

\newtheorem{cor}{Corollary}
\theoremstyle{definition}

\newtheorem{example}{Example}

\newcommand{\emptytab}{\varnothing}

\newcommand{\hcins}{\overset{H}{\longrightarrow}}

\newcommand{\la}{\lambda}

\newcommand{\rank}{\operatorname{rank}}
\newcommand{\word}{\operatorname{word}}
\newcommand{\OO}{{\mathcal O}}
\newcommand{\wb}{\overline}
\newcommand{\bull}{{\scriptscriptstyle \bullet}}
\newcommand{\lft}{{\text{left}}}
\newcommand{\rgt}{{\text{right}}}
\newcommand{\topp}{{\text{top}}}
\newcommand{\bott}{{\text{bottom}}}
\newcommand{\wh}{\widehat}

\newcommand{\pic}[2]{\includegraphics[scale=0.#1]{#2.eps}}

\newcommand{\Groth}{{\mathfrak G}}
\newcommand{\ignore}[1]{}

\psfrag{Ui1i}{$U_{i-1,i}$}
\psfrag{Ai}{$A_i$}
\psfrag{Bi1}{$B_{i-1}$}
\psfrag{Un11}{$U_{n-1,1}$}
\psfrag{Un12}{$U_{n-1,2}$}
\psfrag{Un13}{$U_{n-1,3}$}
\psfrag{U21}{$U_{2,1}$}
\psfrag{U22}{$U_{2,2}$}
\psfrag{U11}{$U_{1,1}$}
\psfrag{Tn}{$T_n$}
\psfrag{T3}{$T_3$}
\psfrag{T2}{$T_2$}
\psfrag{T1}{$T_1$}

\begin{document}

\title[Stable Grothendieck polynomials and $K$-theoretic factor sequences]
{Stable Grothendieck polynomials and $K$-theoretic \\ factor sequences}

\author{Anders Skovsted Buch}
\address{Department of Mathematics, Rutgers University, 110
  Frelinghuysen Road, Piscataway, NJ 08854, USA}
\email{asbuch@math.rutgers.edu}

\author{Andrew Kresch}
\address{Mathematics Institute, University of Warwick,
Coventry CV4 7AL, United Kingdom}
\email{kresch@maths.warwick.ac.uk}

\author{Mark Shimozono}
\address{Department of Mathematics, Virginia Tech, 460 McBryde Hall,
Blacksburg, VA 24061-0123, USA}
\email{mshimo@vt.edu}

\author{\\ Harry~Tamvakis}
\address{Department of Mathematics, Brandeis University - MS 050,
P. O. Box 9110, Waltham, MA 02454-9110, USA}
\email{harryt@brandeis.edu}

\author{Alexander Yong}
\address{Department of Mathematics, University of Minnesota, Minneapolis, 
MN 55455 and The Fields Institute, 222 College Street, Toronto, Ontario,
M5T 3J1, Canada}
\email{ayong@math.umn.edu, ayong@fields.utoronto.ca}

\date{January 19, 2006}


\begin{abstract}
  We formulate a nonrecursive combinatorial rule for the expansion of
  the stable Grothendieck polynomials of [Fomin-Kirillov '94] in the
  basis of stable Grothendieck polynomials for partitions. This gives
  a common generalization, as well as new proofs of the rule of
  [Fomin-Greene '98] for the expansion of the stable Schubert
  polynomials into Schur polynomials, and the $K$-theoretic
  Grassmannian Littlewood-Richardson rule of [Buch '02].  The proof is
  based on a generalization of the Robinson-Schensted and
  Edelman-Greene insertion algorithms. 
  Our results are applied to prove a number of new formulas and
  properties for $K$-theoretic quiver polynomials, and the Grothendieck
  polynomials of [Lascoux-Sch\"{u}tzenberger '82]. In particular, we
  provide the first $K$-theoretic analogue of the factor sequence
  formula of [Buch-Fulton '99] for the cohomological quiver
  polynomials.
\end{abstract}

\maketitle

\markboth{\sc A.~Buch, A.~Kresch, M.~Shimozono, H.~Tamvakis, and A.~Yong}
{\sc Stable Grothendieck polynomials and $K$-theoretic factor sequences}
\pagestyle{myheadings}

\section{Introduction and main results}

\subsection{Stable Grothendieck polynomials}

For each permutation $\pi$ there is a symmetric power series $G_\pi =
G_\pi(x_1,x_2,\dots)$ called the stable Grothendieck polynomial for
$\pi$.  These power series were defined by Fomin and Kirillov
\cite{fomin.kirillov:yang-baxter, fomin.kirillov:grothendieck} as a
limit of the ordinary Grothendieck polynomials of Lascoux and
Sch\"utzenberger \cite{lascoux.schutzenberger:polynomes}.  We recall
this definition in Section \ref{S:groth}.  The term of lowest degree
in $G_\pi$ is the Stanley symmetric function (or stable Schubert
polynomial) $F_\pi$.  The {\em Stanley coefficients\/} which appear in
the Schur expansion of a Stanley function are interesting
combinatorial invariants which generalize the Littlewood-Richardson
coefficients.

Given a partition $\lambda = (\lambda_1 \geq \dots \geq \lambda_k \geq
0)$, the {\em Grassmannian permutation} $\pi_\lambda$ for $\lambda$ is
uniquely defined by the requirement that $\pi_\lambda(i) = i +
\lambda_{k+1-i}$ for $1 \leq i \leq k$ and $\pi_\lambda(i) <
\pi_\lambda(i+1)$ for $i \neq k$.  The power series $G_\lambda :=
G_{\pi_\lambda}$ play a role in combinatorial $K$-theory similar to
the role of Schur functions in cohomology.  Buch has shown
\cite{buch:littlewood-richardson} that any stable Grothendieck
polynomial $G_\pi$ can be written as a finite linear combination
\begin{equation} \label{E:stabcoef}
  G_\pi = \sum_\lambda c_{\pi,\lambda} G_\lambda
\end{equation}
of stable Grothendieck polynomials indexed by partitions, using
integer coefficients $c_{\pi,\lambda}$ that generalize the Stanley
coefficients \cite{buch:grothendieck}.  Lascoux gave a recursive
formula for stable Grothendieck polynomials which confirms a
conjecture that these coefficients have signs that alternate with
degree, i.e.\ $(-1)^{|\lambda|-\ell(\pi)} c_{\pi,\lambda} \geq 0$
\cite{lascoux:transition}.  Here $|\lambda|= \lambda_1+\cdots
+\lambda_k$ and $\ell(\pi)$ is the Coxeter length of $\pi$. The
central result of this paper is a new formula for the coefficients
$c_{\pi,\lambda}$ which generalizes Fomin and Greene's combinatorial
rule \cite{fomin.greene:noncommutative} for Stanley coefficients.

To state our formula, we need the {\em $0$-Hecke monoid}, which is the
quotient of the free monoid of all finite words in the alphabet
$\{1,2,\dotsc\}$ by the relations
\begin{align}
\label{E:idem}
p\,p &\equiv p &\quad&\text{for all $p$} \\
\label{E:braid}
p\,q\,p &\equiv q\,p\,q &\quad&\text{for all $p,q$} \\
\label{E:comm}
p\,q &\equiv q\,p &\quad&\text{for $|p-q|\ge 2$.}
\end{align}
There is a bijection between the $0$-Hecke monoid and the infinite
symmetric group $S_\infty = \bigcup_{n\ge1} S_n$. Given any word $a$
there is a unique permutation $\pi\in S_\infty$ such that $a\equiv b$
for some (equivalently every) reduced word $b$ of $\pi$.  In this case
we write $w(a)=\pi$ and say that $a$ is a {\em Hecke word\/} for
$\pi$.  Notice that the reduced words for $\pi$ are precisely the
Hecke words for $\pi$ that are of minimum length.  Given an additional
permutation $\rho$ with Hecke word $b$, the {\em Hecke product\/} of
$\pi$ and $\rho$ is defined as the permutation $\pi \cdot \rho =
w(ab)$.

We use the English notation for partitions and tableaux.  An
\textit{increasing} tableau is a Young tableau whose rows increase
strictly from left to right, and whose columns increase strictly from
top to bottom. A \textit{decreasing} tableau is defined similarly.
The (column reading) word of a tableau $T$, denoted $\word(T)$, is
obtained by reading the columns of the tableau from bottom to top,
starting with the leftmost column, followed by the column to its
right, etc.  We also define $w(T) := w(\word(T))$.

Our main theorem gives the explicit expansion of the stable
Grothendieck polynomial $G_\pi$ in terms of the $G_\la$.

\begin{thm} \label{T:stabcoef} For any permutation $\pi \in S_\infty$,
  the coefficient $c_{\pi,\lambda}$ in (\ref{E:stabcoef}) is equal to
  $(-1)^{|\lambda|-\ell(\pi)}$ times the number of increasing tableaux
  $T$ of shape $\lambda$ such that $\word(T)$ is a Hecke word for
  $\pi^{-1}$.
\end{thm}

\begin{example} Consider $\pi=31524=s_2s_1s_4s_3$, where each $s_i$ is
  a simple transposition.  The increasing tableaux that provide Hecke
  words for $\pi^{-1}$ are:
\begin{equation*}
\tableau{{1}&{2}\\{3}&{4}}
\qquad \tableau{{1}&{2}&{4}\\{3}}
\qquad \tableau{{1}&{2}&{4}\\{3}&{4}}
\end{equation*}
Hence $G_\pi = G_{22} + G_{31} - G_{32}$.
\end{example}

Theorem \ref{T:stabcoef} may be used to give self-contained proofs of
a number of known results. For example, the finiteness of the
expansion \eqref{E:stabcoef} proved in
\cite{buch:littlewood-richardson} follows immediately from Theorem
\ref{T:stabcoef}.  When the permutation $\pi$ is 321-avoiding,
Theorem~\ref{T:stabcoef} furthermore generalizes Buch's rule for the
coefficients $c_{\pi,\lambda}$ in terms of set-valued tableaux
\cite{buch:littlewood-richardson}, in the sense that there is an
explicit bijection between the relevant increasing and set-valued
tableaux.  As a consequence, we obtain a new proof of the set-valued
Littlewood-Richardson rule for the Schubert structure constants in the
$K$-theory of Grassmannians, as well as an alternative rule based on
increasing tableaux.  This is explained in Section~\ref{S:sv_lr}.

\subsection{Hecke insertion}

Fomin and Kirillov proved that the monomial coefficients of (stable)
Grothendieck polynomials are counted by combinatorial objects called
compatible pairs (also known as resolved wiring diagrams, FK-graphs,
pipe dreams, or nonreduced RC-graphs)
\cite{fomin.kirillov:yang-baxter, fomin.kirillov:grothendieck}.  This
formula was used in \cite{buch:littlewood-richardson} to express the
monomial coefficients of stable Grothendieck polynomials for
partitions in terms of set-valued tableaux (see equation
\eqref{E:Gset}). We prove Theorem~\ref{T:stabcoef} by exhibiting an
explicit bijection between the set of compatible pairs for a
permutation $\pi$ and the set of pairs $(T,U)$ where $T$ is an
increasing tableau with $w(T)=\pi^{-1}$ and $U$ is a set-valued
tableau of the same shape as $T$.  This bijection is constructed using
a new combinatorial algorithm called {\em Hecke insertion}, which is
the technical core of our paper.

Hecke insertion is a generalization of the Edelman-Greene insertion
algorithm \cite{edelman.greene:balanced} (also known as Coxeter-Knuth
insertion) from the set of reduced words to the set of all (Hecke)
words.  It specializes to Robinson-Schensted insertion for words of
distinct integers \cite{robinson, schensted}.  There are two main
novelties in our extension.  First, we need an operation that
``jumps'' many columns at once.  Second, an accompanying reverse
insertion algorithm can pass back different intermediate values than
the insertion algorithm generated.  Neither of these elements appear
in the classical algorithms.  We also use Hecke insertion to define
products of decreasing tableaux, which enter into our definition of
$K$-theoretic factor sequences.

\subsection{Quiver varieties}
\label{S:intro_quiver}

Our main application of Theorem~\ref{T:stabcoef} concerns the classes
of quiver varieties in $K$-theory. Recall that a sequence of vector
bundle morphisms $E_0 \to E_1 \to \dots \to E_n$ over a non-singular
variety ${\mathfrak X}$ together with a set of {\em rank conditions\/}
$r=\{r_{ij}\}$ for $0\leq i \leq j \leq n$ define a {\em quiver
  variety\/} $\Omega_r \subset {\mathfrak X}$ of points where each
composition of bundle maps $E_i \to E_j$ has rank at most $r_{ij}$.
We demand that the rank conditions can actually {\em occur}, and that
the bundle maps are generic, so that the quiver variety $\Omega_r$
obtains its expected codimension $d(r) = \sum_{i<j} (r_{i,j-1} -
r_{ij})(r_{i+1,j}-r_{ij})$.  Buch and Fulton proved a formula for the
cohomology class of $\Omega_r$ \cite{buch.fulton:chern}, which was
later generalized to $K$-theory by Buch \cite{buch:grothendieck}.  The
$K$-theory version states that the Grothendieck class of $\Omega_r$ is
given by
\begin{equation} \label{E:kquiver}
  [\OO_{\Omega_r}] = \sum_\mu c_\mu(r)\,
  G_{\mu_1}(E_1-E_0) G_{\mu_2}(E_2-E_1) \cdots G_{\mu_n}(E_n-E_{n-1}) \,,
\end{equation}
where the sum is over sequences $\mu = (\mu_1,\dots,\mu_n)$ of
partitions $\mu_i$ such that $\sum |\mu_i| \geq d(r)$ and each
partition $\mu_i$ can be contained in the rectangle $e_i \times
e_{i-1}$ with $e_i$ rows and $e_{i-1}$ columns, where $e_i := r_{ii} =
\rank(E_i)$. The notation $G_{\mu_i}(E_{i+1}-E_i)$ will be explained
in Section~2.

The coefficients $c_\mu(r)$ in formula (\ref{E:kquiver}) are integers
called {\em quiver coefficients}.  When $\sum |\mu_i| = d(r)$, the
coefficient $c_\mu(r)$ also appears in the cohomology formula from
\cite{buch.fulton:chern} and is called a {\em cohomological quiver
  coefficient}. A precise conjecture for these cohomological
coefficients was posed in \cite{buch.fulton:chern}, which asserts that
$c_\mu(r)$ counts the number of {\em factor sequences\/} of tableaux
with shapes given by the sequence of partitions $\mu$. A factor
sequence is a sequence of semistandard Young tableaux that can be
obtained by performing a series of plactic factorizations and
multiplications of chosen tableaux arranged in a {\em tableau
  diagram}.  For a specific choice of tableau diagram, this conjecture
was proved by Knutson, Miller and Shimozono
\cite{knutson.miller.ea:four}.  However, the original definition of
factor sequences from \cite{buch.fulton:chern} as sequences of
tableaux generated using the plactic product has no known
generalization to $K$-theory.

In this paper, we prove that $K$-theoretic quiver coefficients are
counted by a new type of factor sequence, generalizing the
cohomological factor sequences defined by Buch in
\cite{buch:alternating} using the Coxeter-Knuth product of tableaux.
The new $K$-theoretic factor sequences are constructed from a tableau
diagram of decreasing\footnote{The combinatorics of these factor
  sequences naturally requires decreasing rather than increasing
  tableaux.}  tableaux using the same algorithm that defines the
original factor sequences, except that the plactic product is replaced
with a product $(U,T) \mapsto U \cdot T$ of decreasing tableaux which
is compatible with Hecke products of permutations (see Section
\ref{SS:product}).

For each $0 \leq i < j \leq n$ let $R_{ij}$ be a rectangle with
$r_{i+1,j}-r_{ij}$ rows and $r_{i,j-1}-r_{ij}$ columns.  Let $U_{ij}$
be the unique decreasing tableau of shape $R_{ij}$ such that the lower
left box contains the number $r_{i,j-1}$, and the number in each box
is one larger than the number below it and one smaller than the number
to the left of it.  For example, if $r_{i,j-1} = 6$, $r_{i+1,j} = 5$,
and $r_{ij} = 2$ then
\[ 
U_{ij} \ = \ \ \tableau{{8}&{7}&{6}&{5}\\{7}&{6}&{5}&{4}\\{6}&{5}&{4}&{3}} \ . 
\]
These tableaux $U_{ij}$ can be arranged in a triangular tableau
diagram as in \cite[\S 4]{buch.fulton:chern}:
\[ \begin{matrix}
U_{01} && U_{12} && \cdots && U_{n-1,n} \\
& U_{02} && \cdots && U_{n-2,n} \\
&& \ddots \\
&&& U_{0n}
\end{matrix} \]

We define a {\em $K$-theoretic factor sequence\/} for the rank
conditions $r$ by induction on $n$.  If $n=1$ then the only factor
sequence is the sequence $(U_{01})$ consisting of the only tableau in
the tableau diagram.  If $n \geq 2$ then the numbers $\wb r = \{ \wb
r_{ij} : 0 \leq i \leq j \leq n-1 \}$ defined by $\wb r_{ij} =
r_{i,j+1}$ form a valid set of rank conditions corresponding to a
sequence of $n-1$ bundle maps.  In this case, a factor sequence for
$r$ is any sequence of the form $(U_{01} \cdot A_1, \dots, B_{i-1}
\cdot U_{i-1,i} \cdot A_i, \cdots, B_{n-1} \cdot U_{n-1,n})$, for a
choice of decreasing tableaux $A_i$ and $B_i$ such that $(A_1\cdot
B_1, \dots, A_{n-1}\cdot B_{n-1})$ is a factor sequence for $\wb r$.

\begin{thm} \label{T:quivcoef}
  The $K$-theoretic quiver coefficient $c_\mu(r)$ is equal to
  $(-1)^{\sum |\mu_i| - d(r)}$ times the number of $K$-theoretic
  factor sequences $(T_1,\dots,T_n)$ for the rank conditions $r$, such
  that $T_i$ has shape $\mu_i$ for each $i$.
\end{thm}

Central to the proof of the nonnegativity of cohomological quiver
coefficients given in \cite{knutson.miller.ea:four} is the {\em stable
  component formula}, which writes the cohomology class of a quiver
variety as a sum of products of Stanley functions.  This sum is over
all {\em lace diagrams\/} representing the rank conditions $r$, which
have the smallest possible number of crossings.  The $K$-theoretic
version of the component formula from \cite{buch:alternating,
  miller:alternating} states that
\begin{equation} \label{E:stabcomp}
  [\OO_{\Omega_r}] = \sum_{(\pi_1,\dots,\pi_n)}
  (-1)^{\sum \ell(\pi_i) - d(r)}
  G_{\pi_1}(E_1-E_0) \cdots G_{\pi_n}(E_n-E_{n-1})
\end{equation}
where the sum is over a generalization of minimal lace diagrams, which
was named {\em KMS-factorizations\/} in \cite{buch:alternating}.  We
recall this definition in Section \ref{S:quiver}.  Our proof of
Theorem~\ref{T:quivcoef} is based on Theorem~\ref{T:stabcoef},
equation (\ref{E:stabcomp}), and the following characterization of the
$K$-theoretic factor sequences:

\begin{thm} \label{T:fackms}
  A sequence of decreasing tableaux $(T_1,\dots,T_n)$ is a
  $K$-theoretic factor sequence for the rank conditions $r$ if and
  only if $(w(T_1),\dots,w(T_n))$ is a KMS-factorization for $r$.
\end{thm}

\subsection{Outline of the paper}

We give a brief overview of the contents of this article.  In Section
\ref{S:groth} we recall the original definitions of Grothendieck
polynomials and stable Grothendieck polynomials, and record two useful
monomial expansions for the latter.  In Section \ref{S:Hecke} we
define the Hecke insertion algorithm, establish its basic properties,
give the proof of Theorem \ref{T:stabcoef}, and apply it to reprove
the set-valued Littlewood-Richardson rule from
\cite{buch:littlewood-richardson}.  We also give an algorithm for
generating all increasing tableaux which represent a given
permutation.  Finally, our applications to quiver coefficients and
factor sequences are contained in Section \ref{S:quiver}. These
include a proof that the $K$-theoretic quiver coefficients are special
cases of the coefficients $c_{z(r),\lambda}$ in the expansion
\eqref{E:stabcoef} of the stable Grothendieck polynomial for a
Zelevinsky permutation $z(r)$ (Theorem \ref{T:quivstab}), and a new
formula for the decomposition coefficients of universal Grothendieck
polynomials (Theorem \ref{T:univgroth}).

\section{Grothendieck polynomials} \label{S:groth}

Grothendieck polynomials were introduced by Lascoux and
Sch\"utzenberger \cite{lascoux.schutzenberger:structure} as polynomial
representatives for the classes of structure sheaves of Schubert
varieties in the $K$-theory of the flag variety for $GL_n$.  Let
$X=(x_1,x_2,\ldots)$ and $Y=(y_1,y_2,\ldots)$ be two sequences of
commuting independent variables and $\pi \in S_{n}$. If $\pi = \pi_0$
is the longest permutation in $S_n$, then we set
\[ \Groth_{\pi_0}(X;Y) = \prod_{i+j\leq n} (x_i+y_j-x_i y_j) \,. \] If
$\pi\neq \pi_0$, we can find a simple transposition $s_i=(i, i+1)\in
S_{n}$ such that $\ell(\pi s_i)=\ell(\pi)+1$.  We then define
\[ \Groth_\pi(X;Y) = \frac{(1-x_{i+1})\Groth_{\pi
    s_i}(X ; Y) - (1-x_i)\Groth_{\pi
    s_i}(x_1,\ldots,x_{i+1},x_i,\ldots,x_n; Y)}{x_i-x_{i+1}} \,.
\]

For $\pi\in S_\infty$ and $r\ge0$, let $1^r \times \pi\in S_\infty$
denote the permutation obtained by putting $r$ fixed points in front
of $\pi$, that is, $1^r \times \pi = \rho$ where $\rho(i)=i$ for $1\le
i\le r$ and $\rho(i)=\pi(i-r)+r$ for $i>r$. The {\em stable
  Grothendieck polynomial} $G_\pi(X;Y)$ is the formal power series,
symmetric in the $X$ and $Y$ variables separately, defined by
\[ G_\pi(X;Y) = \lim_{r\to\infty} \Groth_{1^r \times \pi}(X;Y) \,. \]
Given vector bundles $E = L_1\oplus\cdots\oplus L_p$ and $F =
M_1\oplus\cdots \oplus M_q$ over a variety ${\mathfrak X}$ which are
direct sums of line bundles, we write
\[ G_{\lambda}(E-F)=G_{\lambda}(1-L_1^{-1},\ldots, 1-L_p^{-1}; 
   1-M_1,\ldots,1-M_q) \in K({\mathfrak X})\, .
\]
The symmetry of $G_{\lambda}(X;Y)$ implies that this is a polynomial
in the exterior powers of $E^\vee$ and $F$.  Therefore
$G_{\lambda}(E-F)$ makes sense even for bundles that are not split
into direct sums of line bundles. This explains the notation used in
(\ref{E:kquiver}).

We will be mostly interested in the specialization $G_\pi =
G_{\pi}(X;0)$.  We recall Fomin and Kirillov's combinatorial
construction of these polynomials \cite{fomin.kirillov:grothendieck},
using notation which generalizes Billey, Jockusch, and Stanley's
formula for Schubert polynomials \cite{billey.jockusch.ea:schub}.
Define a {\em compatible pair\/} to be a pair $(a,i)$ of words $a=a_1
a_2 \cdots a_p$ and $i = i_1 i_2 \cdots i_p$ of positive integers,
such that $i_1 \leq i_2 \leq \dots \leq i_p$, and so that
$a_j>a_{j+1}$ whenever $i_j=i_{j+1}$.  The stable Grothendieck
polynomial for $\pi \in S_\infty$ is then given by
\cite{fomin.kirillov:grothendieck}
\begin{equation} \label{E:Groth}
  G_\pi = \sum_{(a,i)} (-1)^{\ell(i)-\ell(\pi)}\, x^i
\end{equation}
where the sum is over all compatible pairs $(a,i)$ such that
$w(a)=\pi$.  Here $\ell(i)$ is the common length of $a$ and $i$, and
$x^i = x_{i_1} x_{i_2} \cdots x_{i_{\ell(i)}}$.

\ignore{
Say that a compatible pair is \textit{bounded} if $i_j \le a_j$ for
all $j$. The Grothendieck polynomial $\Groth_\pi(X)$ is given by the
sum \eqref{E:Groth} except that $(a,i)$ must also be bounded. The
\textit{Schubert polynomial} for $\pi$ is equal to the lowest degree
term of $\Groth_\pi$, while the \textit{Stanley function} or
\textit{stable Schubert polynomial} $F_\pi$ is the symmetric function
given by the lowest degree term of $G_\pi$.
}

A \textit{set-valued tableau} of shape $\la$ is a filling of the boxes
of the Young diagram of $\la$ with finite nonempty sets of positive
integers, such that these sets are weakly increasing along rows and
strictly increasing down columns.  In other words, all integers in a
box must be smaller than or equal to the integers in the box to the
right of it, and strictly smaller than the integers in the box below
it.  For a set-valued tableau $S$, let $x^S$ denote the monomial where
the exponent of $x_i$ is equal to the number of boxes containing the
integer $i$, and let $|S|$ be the degree of this monomial. Buch's
formula for the monomial expansion of $G_\la$ is given by
\cite{buch:littlewood-richardson}
\begin{equation} \label{E:Gset}
  G_\la = \sum_S (-1)^{|S|-|\lambda|}\, x^S
\end{equation}
where $S$ runs over all set-valued tableaux of shape $\la$.

\section{Hecke Insertion and the proof of Theorem~\ref{T:stabcoef}}
\label{S:Hecke}

In view of \eqref{E:Groth} and \eqref{E:Gset}, to prove
Theorem~\ref{T:stabcoef} it suffices to establish a bijection $(a,i)
\mapsto (T,U)$ between all compatible pairs $(a,i)$ such that $w(a) =
\pi$, and all pairs of tableaux $(T,U)$ of the same shape, such that
$T$ is increasing with $w(T) = \pi^{-1}$ and $U$ is set-valued.  In
addition, this bijection must satisfy $x^U = x^i$.  To construct this
bijection, we need a new algorithm called {\em Hecke insertion}.

\subsection{Hecke insertion} \label{SS:heckeins}

We shall define the Hecke (column) insertion of a non-negative integer
$x$ into the increasing tableau $Y$, resulting in the increasing
tableau $Z$.  The shape of $Z$ always contains the shape of $Y$, and
contains at most one extra box $c$.  Unlike ordinary
Robinson-Schensted insertion, it is possible that $Z$ has the same
shape as $Y$, but even in this case it will contain a special corner
$c$ where the insertion algorithm terminated.  To keep track of these
cases, we will use a parameter $\alpha \in \{0,1\}$, which is set to 1
if and only if the corner $c$ is outside the shape of $Y$.  Thus the
complete output of the insertion algorithm is the triple
$(Z,c,\alpha)$.  We will use the notation $Z = (x \hcins Y)$.

The algorithm proceeds by inserting the integer $x$ into the first
column of $Y$.  This may modify this column, and possibly produce an
{\em output integer}, which is then inserted into the second column of
$Y$, etc.  This process is repeated until an insertion does not
produce an output integer.  The procedure for inserting an integer $x$
into a column $C$ is as follows.

If $x$ is larger than or equal to all boxes of $C$, then no new output
value is produced and the algorithm terminates.  If adjoining $x$ as a
new box below $C$ results in an increasing tableau, then $Z$ is the
resulting tableau, $\alpha=1$, and $c$ is the new corner where $x$ was
added.  If $x$ cannot be added, then no further modifications are
carried out to produce $Z$, $\alpha=0$, and $c$ is the corner of the
row of $Z$ containing the bottom box of the column $C$.

Otherwise $C$ contains boxes strictly larger than $x$, and we let $y$
be the smallest such box.  If replacing $y$ with $x$ results in an
increasing tableau, then this is done.  In either case, $y$ is the
output integer, which is inserted into the next column.

\begin{example} \label{X:1}
\begin{equation*}
  3 \hcins \ \tableau{{1}&{2}&{3}&{4}\\{2}&{5}&{6}\\{3}\\{5}}
  \ = \ \tableau{{1}&{2}&{3}&{4}\\{2}&{5}&{6}\\{3}\\{5}}
\end{equation*}
The integer $3$ is inserted into the first column, which contains $3$.
So $5$ is inserted into the second column, whose largest value is $5$.
The algorithm terminates with $\alpha=0$, and $c=(2,3)$ is the corner
in the second row and third column.
\end{example}

\begin{example} \label{X:2a}
\begin{equation*}
  2 \hcins \ \tableau{{1}&{2}\\{2}&{5}\\{4}}
  \ = \ \tableau{{1}&{2}&{5}\\{2}&{4}\\{4}}
\end{equation*}
The integer $2$ is inserted into the first column, which contains $2$.
So $4$ is inserted into the second column, displacing the $5$.  The
$5$ is inserted into the third column, where it comes to rest.  We get
$\alpha=1$ and $c=(1,3)$.
\end{example}

\begin{example} \label{X:2b}
\begin{equation*}
  2 \hcins \ \tableau{{1}&{2}&{3}&{5}\\{2}&{3}&{4}\\{4}}
  \ = \ \tableau{{1}&{2}&{3}&{5}\\{2}&{3}&{4}\\{4}}
\end{equation*}
The integer $2$ is inserted into the first column, which contains a
$2$.  So $4$ is inserted into the second column, which has largest
entry $3$.  Since the first column still contains the value $4$ in its
bottom box, it is not possible to add a box with $4$ to the second
column.  We obtain $\alpha=0$, and $c = (2,3)$ is the corner of the
second row.
\end{example}

\begin{example} \label{X:3}
\begin{equation*}
  1 \hcins \ \tableau{{1}&{2}&{3}\\{3}&{4}&{5}}
  \ = \ \tableau{{1}&{2}&{3}&{5}\\{3}&{4}&{5}}
\end{equation*}
The integer $1$ is inserted into the first column, which already
contains a $1$.  So $3$ is inserted into the second column.  It would
have replaced $4$, but this replacement would place a $3$ directly to
the right of another $3$, violating the increasing tableau condition.
So the second column is unchanged and $4$ is inserted into the third
column.  Similarly $4$ cannot replace $5$, so $5$ is inserted into the
fourth column, where it comes to rest in the cell $c=(1,4)$ with
$\alpha=1$.
\end{example}

\subsection{Reverse Hecke insertion}

Let $Z$ be an increasing tableau, $c$ a corner of $Z$, and $\alpha \in
\{0,1\}$.  Reverse Hecke insertion applied to the triple
$(Z,c,\alpha)$ produces a pair $(Y,x)$ of an increasing tableau $Y$
and a positive integer $x$ as follows.  Let $y$ be the integer in the
cell $c$ of $Z$.  If $\alpha=1$ then remove $y$.  In any case, reverse
insert $y$ into the column to the left of the corner $c$.

Whenever a value $y$ is reverse inserted into a column $C$, let $x$ be
the largest entry of $C$ such that $x < y$.  If replacing $x$ with $y$
results in an increasing tableau, then this is done.  In any case, the
integer $x$ is passed along to the left.  If $C$ is not the left-most
column, this means that $x$ is reverse inserted into the column left
of $C$; otherwise $x$ becomes the final output value, along with the
modified tableau.

\begin{example} \label{X:revins} Let us apply reverse Hecke insertion
  to the tableau computed in Example \ref{X:3} at the cell $c=(1,4)$
  with $\alpha=1$.  The integer $5$ in this cell is then removed, and
  $5$ is reverse inserted into the third column.  Since $5$ is already
  in the third column, it is not changed, and $3$ is reverse inserted
  into the second column.  Here $3$ cannot replace $2$ because this
  would place a $3$ directly to the left of a $3$.  The second column
  is unchanged and $2$ is reverse inserted into the first column.  The
  $2$ cannot replace $1$ for the same reason, so the first column is
  unchanged and $x = 1$ is the output value.  This recovers the
  initial tableau of Example \ref{X:3}.
\end{example}

Let ${\mathcal I}$ denote the set of pairs $(Y,x)$ where $Y$ is an
increasing tableau and $x$ is a positive integer.  Let ${\mathcal R}$
be the set of triples $(Z,c,\alpha)$ where $Z$ is an increasing
tableau, $c$ a corner cell of $Z$, and $\alpha \in \{0,1\}$.

\begin{thm} 
\label{T:bij} 
Hecke insertion $(Y,x) \mapsto (Z,c,\alpha)$ and reverse Hecke
insertion $(Z,c,\alpha) \linebreak \mapsto (Y,x)$ define mutually
inverse bijections between the sets ${\mathcal I}$ and ${\mathcal R}$.
\end{thm}
\begin{proof}
  Assume at first that $(Z,c,\alpha)$ has been obtained by applying
  Hecke insertion to $(Y,x)$.  We must show that reverse Hecke
  insertion recovers $(Y,x)$ from $(Z,c,\alpha)$.  This is clear if
  $Y$ is the empty tableau.  In general we proceed by induction on the
  number of columns in $Y$.

  If $x$ is strictly larger than all the integers in the first column
  of $Y$, then $Z$ is obtained by adding a box containing $x$ to the
  first column of $Y$, $c$ is this box, $\alpha = 1$, and reverse
  Hecke insertion clearly maps $(Z,c,\alpha)$ back to $(Y,x)$.  If $x$
  is equal to the largest integer in the first column of $Y$, then
  $Z=Y$, $\alpha=0$, and $c$ is the leftmost corner of $Y$.  Also in
  this case it is easy to see that reverse Hecke insertion applied to
  $(Z,c,\alpha)$ recovers $(Y,x)$.

  We can therefore assume that the first column of $Y$ contains at
  least one integer that is strictly larger than $x$.  We let $Y'$
  denote the first column of $Y$ and let $Y''$ be the rest of $Y$.  We
  define $Z'$ and $Z''$ similarly, and we regard $c$ as a corner of
  both $Z$ and of $Z''$.  Let $a$ be the smallest box of $Y'$ for
  which $x < a$.  If $(Z'',c,\alpha)$ is the result of applying Hecke
  insertion to $(Y'',a)$, then we know by induction that reverse Hecke
  insertion applied to $(Z,c,\alpha)$ first recovers $Y''$ from $Z''$,
  after which $a$ is reverse inserted into $Z'$.  If the box of $Y'$
  containing $a$ was replaced with $x$ in the construction of $Z'$,
  then the reverse insertion puts $a$ back in this box to recover
  $Y'$, and the output value of the reverse insertion is $x$.  If the
  box of $Y'$ containing $a$ was not changed to $x$ when $Z'$ was
  formed, then $Z'=Y'$ must contain $x$ in the box immediately over
  $a$.  Reverse inserting $a$ therefore leaves $Z'$ unchanged, and
  gives $x$ as the output value, as required.

  It remains to consider the case when $Z''$ differs from the tableau
  $(a \hcins Y'')$.  This can only happen if the box of $Y'$
  containing $a$ was not replaced with $x$, and if the leftmost box in
  the same row of $(a \hcins Y'')$ contains $a$.  Let $b$ be the
  largest box of the first column of $Y''$ for which $b \leq a$, and
  observe that $(Z'',c,\alpha)$ must be the result of applying Hecke
  insertion to $(b,Y'')$.  It therefore follows by induction that
  reverse Hecke insertion applied to $(Z,c,\alpha)$ first recovers
  $Y''$ from $Z''$, after which $b$ is reverse inserted into $Z'=Y'$.
  Finally, notice that $Y'$ must contain $x$ in the box above the box
  containing $a$, since otherwise $a$ would have been replaced with
  $x$ in the initial Hecke insertion.  Furthermore we have $x < b$,
  since $b$ is in the box of $Y''$ to the right of $x$ in $Y'$.  Since
  $x < b \leq a$, we conclude that reverse insertion of $b$ to $Y'$
  leaves this column unchanged, and $x$ is the final output value, as
  required.  This verifies that Hecke insertion followed by reverse
  Hecke insertion is the identity map.

  To finish the proof, assume that $(Y,x)$ has been obtained by
  applying reverse Hecke insertion to $(Z,c,\alpha)$.  We must show
  that Hecke insertion maps $(Y,x)$ back to $(Z,c,\alpha)$.  This is
  easily checked if $c$ is the corner of the bottom row of $Z$.  In
  general we use induction on the number of columns in $Z$.  Assume
  that $c$ is not in the bottom row of $Z$, and let $Y'$, $Y''$, $Z'$,
  and $Z''$ be defined as above.  We know that for some integer $a$,
  the pair $(Y'',a)$ is the result of applying reverse Hecke insertion
  to $(Z'',c,\alpha)$.  We also know that $Y'$ is obtained by reverse
  inserting $a$ into $Z'$, keeping in mind that $Y''$ resides to the
  right, and $x$ is the output value resulting from this insertion.

  We first prove that $Y'$ contains integers strictly larger than $x$.
  In fact, if this was not true, then $x$ would be equal to the bottom
  box of $Y'=Z'$.  Since $a$ did not replace $x$ in $Z'$, we deduce
  that $Z'$ and $Y''$ have equally many rows, and $a$ is the content
  of the bottom-left box of $Y''$.  On the other hand, since $c$ is
  not in the bottom row of $Z''$, it follows by induction on the
  number of columns of $Z$ that the first column of $Y''$ contains
  integers strictly larger than $a$, a contradiction.

  If the box of $Z'$ containing $x$ is replaced by $a$ in $Y'$, then
  Hecke insertion applied to $(Y,x)$ first restores $Z'$ in the first
  column, after which $a$ is inserted into the first column of $Y''$.
  Since we know by induction that $(Z'',c,\alpha)$ is the result of
  applying Hecke insertion to $(a,Y'')$, we deduce that Hecke
  insertion maps $(Y,a)$ to $(Z,c,\alpha)$, as required.

  Finally assume that $Y'=Z'$, i.e.\ $x$ is not replaced by $a$.
  Since $Y'$ contains integers strictly larger than $x$, we know that
  $x$ is not in the bottom box of $Y'$.  Let $y$ be the integer in the
  box of $Y'$ just under $x$.  When $x$ is Hecke inserted into
  $Y'=Z'$, this column is not changed, and $y$ is inserted into $Y''$.
  If $y=a$, then we know by induction that this recovers
  $(Z'',c,\alpha)$.  Otherwise we must have $y > a$, and $a$ must be
  contained in the box in the first column of $Y''$ which is in the
  same row as the box containing $x$ in $Y'$.  In this case Hecke
  insertion of $y$ into $Y''$ with $Z'$ to the left will produce the
  same result as Hecke inserting $a$ into $Y''$ with nothing to the
  left, namely $(Z'',c,\alpha)$.  This verifies that reverse Hecke
  insertion followed by Hecke insertion is the identity map, which
  completes the proof.
\end{proof}

\subsection{Properties of Hecke insertion}

We will need two additional properties of Hecke insertion.  The first
property says that Hecke insertion respects Hecke words.

\begin{lemma} \label{L:Hecke}
  Let $Y$ be an increasing tableau, $x$ a positive integer, and set $Z
  = (x \hcins Y)$.  Then $\word(Z) \equiv x\, \word(Y)$.
\end{lemma}
\begin{proof}
  It is easiest to check that reverse Hecke insertion preserves Hecke
  words.  It is enough to consider the following situation.  Let $C$
  be a column, $y$ a number that is reverse inserted into $C$, $U$ the
  modified tableau to the right of $C$ from which $y$ comes, $C'$ the
  modification of $C$, and $x$ the output value. We must show that
  \[ \word(C)\, y\, \word(U) \equiv x\, \word(C')\, \word(U) \,. \]

  If $x$ is replaced by $y$ in $C$ then we have $\word(C)\, y \equiv
  x\, \word(C')$ by the relations \eqref{E:comm}.  If $x$ is not
  replaced because $y$ is immediately below $x$ in $C$, then
  $\word(C)\, y \equiv x\, \word(C')$ holds by the relations
  \eqref{E:braid} and \eqref{E:comm}.

  Finally, assume that $x$ is not replaced by $y$ in $C$ because the
  box of $C$ containing $x$ is just left of a box in $U$ containing
  $y$.  In this case we show that $\word(C) = \word(C') \equiv x\,
  \word(C')$ and $y\, \word(U) \equiv \word(U)$.  If $x$ is contained
  in the bottom box of $C$, then the first relation follows from
  (\ref{E:idem}).  Otherwise let $z$ be the box just under $x$ in $C$.
  Since we must have $x < y < z$, the relation $\word(C') \equiv x\,
  \word(C')$ follows from (\ref{E:idem}) and (\ref{E:comm}).
  Similarly, if the first column of $U$ contains an integer $w$ in a
  box just below $y$, then $y < z < w$, and (\ref{E:idem}) and
  (\ref{E:comm}) imply that $y\, \word(U) \equiv \word(U)$.
\end{proof}

We also need the following ``Pieri property'' of Hecke insertion.

\begin{lemma} \label{L:pieri} Let $Y$ be an increasing tableau, and
  $x_1, x_2$ two positive integers.  Suppose that Hecke insertion of
  $x_1$ into $Y$ results in $(Z,c_1,\alpha_1)$ and that Hecke
  insertion of $x_2$ into $Z$ results in $(T,c_2,\alpha_2)$.  Then
  $c_2$ is strictly to the right of $c_1$ if and only if $x_1 > x_2$.
\end{lemma}
\begin{proof}
  As in the proof of Lemma~\ref{L:Hecke}, it is easier to work with
  reverse Hecke insertion.  We consider $x_2$ as the output value
  obtained from applying reverse Hecke insertion to $T$ starting at
  the corner $c_2$, and $x_1$ as the output value obtained by applying
  reverse insertion to the corner $c_1$ of the result.  We first
  consider the case that $c_2$ is strictly to the right of $c_1$.
  
  Suppose $c_1$ is in the first column of $T$.  If the path of the
  first reverse insertion went through $c_1$, then the number in $c_1$
  became larger, so the lemma holds.  If the first insertion path went
  above $c_1$, then clearly the lemma also holds.
  
  Consider $T$ as split vertically into the subtableau of columns
  weakly to the right of $c_1$ and the subtableau of columns strictly
  to the left. The above observations imply that it is enough to prove
  the following: Let $y_2 < y_1$ and assume $y_2$ is first reverse
  inserted into a column $C$, and then $y_1$ is reverse inserted into
  the modification of $C$, then the first output value is strictly
  smaller than the second.
  
  The first output value $x$ is the bottom-most entry in $C$ which is
  strictly smaller than $y_2$.  If the first reverse insertion
  replaces $x$ by $y_2$ or if the box below $x$ contains $y_2$, then
  the second reverse insertion produces an output value which is
  greater than or equal to $y_2$.
  
  Finally assume that when the first reverse insertion occurs, $y_2$
  is contained in a box just to the right of the box of $C$ containing
  $x$.  In this case $y_1$ must reside in a box below $y_2$ in the
  column to the right of $C$, which implies that the second output
  value will come from a box below $x$ in $C$, and thus be strictly
  larger than $x$.

  Now we consider the case when $c_2$ is further to the left or in the
  same column as $c_1$.  We must show that the first output value is
  larger than or equal to the second. If $c_2$ is in the first column
  then this is clear. The first output value is the largest entry in
  the first column.  After the reverse insertion from $c_2$, the
  second output value must come from the first column, all of whose
  entries are smaller than or equal to the first output value.
  
  For reasons similar to those of the previous case, we need to show
  the following: suppose $y_2 \geq y_1$ and that $y_2$ is reverse
  inserted into a column $C$, and then $y_1$ is reverse inserted into
  the modification of $C$, then the first output value is larger than
  or equal to the second.
  
  The first output value $x$ is the bottom-most entry of $C$ that is
  strictly smaller than $y_2$. If $x$ is replaced by $y_2$, then the
  bottom-most entry of the modified column which is strictly smaller
  than $y_1$ must be located above the current location of $y_2$.  And
  this entry is still smaller than $x$.
  
  Suppose $x$ is not replaced by $y_2$. In this case, the bottom-most
  entry that is strictly smaller than $y_1$ is either $x$ or something
  above $x$. In both cases, the second output value is less than or
  equal to $x$.
\end{proof}

\subsection{Proof of Theorem~\ref{T:stabcoef}} 
\label{SS:Column} 

Recall from the discussion at the beginning of Section \ref{S:Hecke}
that to prove Theorem \ref{T:stabcoef}, it suffices to exhibit a
bijection $(a,i)\mapsto (T,U)$ where $(a,i)$ is a compatible pair with
$w(a)=\pi$, $T$ is an increasing tableau with $w(T)=\pi^{-1}$, and $U$
is a set-valued tableau of the same shape as $T$.  Moreover this
bijection must satisfy $x^i=x^U$.

Let $(a,i)$ be as above with $a=a_1\dotsm a_p$ and $i=i_1\dotsm i_p$.
We start with the empty tableau pair $(T_0,U_0) = (\emptytab,
\emptytab)$.  If $(T_{j-1},U_{j-1})$ has been defined for some $j\geq
1$, let $(T_j,c_j,\alpha_j)$ be the result of Hecke inserting $a_j$
into $T_{j-1}$.  If $\alpha_j = 1$ then $U_j$ is obtained from
$U_{j-1}$ by adding the corner $c_j$ and putting the singleton set
$\{i_j\}$ in this box.  Otherwise $c_j$ is already a corner of
$U_{j-1}$, and $U_j$ is obtained by putting $i_j$ into the existing
set in this corner of $U_{j-1}$.  We finally set $(T,U) = (T_p,U_p)$.

The map $(a,i)\mapsto (T,U)$ has the desired properties.  $U$ is a
set-valued tableau by Lemma~\ref{L:pieri} and $x^i=x^U$ by definition.
The fact that $w(T)=\pi^{-1}$ follows from Lemma \ref{L:Hecke},
combined with the fact that the reversal of words gives a bijection
between the Hecke words for $\pi$ and those for $\pi^{-1}$.

Finally, for $j \geq 1$ we note that $i_j$ is the largest integer
appearing in $U_j$, and $c_j$ is the (unique) rightmost corner of
$U_j$ containing this integer.  If this corner of $U_j$ contains only
a singleton, then $\alpha_j=1$, and $U_{j-1}$ can be obtained by
removing the box $c_j$ from $U_j$.  Otherwise we have $\alpha_j=0$ and
$U_{j-1}$ is obtained by removing the integer $i_j$ from the box $c_j$
of $U_j$.  Since the pair $(T_{j-1},a_j)$ is the result of applying
reverse Hecke insertion to the triple $(T_j,c_j,\alpha_j)$, this shows
that the integers $i_j$ and $a_j$ and the pair $(T_{j-1},U_{j-1})$ can
be recovered from $(T_j,U_j)$.  Repetition of this procedure provides
an inverse map $(T,U) \mapsto (a,i)$.  This completes the proof of
Theorem \ref{T:stabcoef}.

\begin{example} \label{E:compat} Let $(a,i)$ be the compatible pair
  given by $a=41443$ and $i=11244$.  Then the above proof constructs
  the following sequence of tableau pairs $(T_j,U_j)$.
\begin{equation*}
\left(\,\tableau{{4}}\,,\,\tableau{{1}}\,\right),
\left(\,\tableau{{1}&{4}}\,,\,\tableau{{1}&{1}}\,\right),
\left(\,\tableau{{1}&{4}\\{4}}\,,\,\tableau{{1}&{1}\\{2}}\,\right),
\left(\,\tableau{{1}&{4}\\{4}}\,,\,\tableau{{1}&{1}\\{24}}\,\right),
\left(\,\tableau{{1}&{4}\\{3}}\,,\,\tableau{{1}&{14}\\{24}}\,\right)
\end{equation*}
\end{example}

\subsection{Increasing tableaux and set-valued tableaux} \label{S:sv_lr}

In this section we sketch how to recover the set-valued
Littlewood-Richardson rule of \cite{buch:littlewood-richardson} from
Theorem~\ref{T:stabcoef}.

A compatible pair $(a,i)$ can be identified with a diagram of columns
of boxes containing positive integers, which increase strictly from
top to bottom.  The boxes of column $p$ contain the integers $a_j$ for
which $i_j = p$.  Some of the columns may be empty.  For example, the
compatible pair from Example~\ref{E:compat} is represented by the
diagram:
\[
\begin{picture}(0,0)
\put(0,11){\line(1,0){44}}
\end{picture}
\tabl{{1}&{4}&&{3}\\{4}&&&{4}}
\]

Let $\lambda/\mu$ be the skew Young diagram between two partitions
$\lambda$ and $\mu$.  If we let $\pi_\lambda$ and $\pi_\mu$ be the
corresponding Grassmannian permutations with descent at the same
position $k$, then the 321-avoiding permutation corresponding to the
skew shape $\lambda/\mu$ is given by $\pi_{\lambda/\mu} =
\pi_\lambda\,\pi_\mu^{-1}$ (see \cite{billey.jockusch.ea:schub}).  It
was shown in \cite[Thm.~3.1]{buch:littlewood-richardson} that the
stable Grothendieck polynomial $G_{\lambda/\mu} :=
G_{\pi_{\lambda/\mu}}$ can be written as the sum
\[ G_{\lambda/\mu} = \sum_S (-1)^{|S|-|\lambda/\mu|} x^S \]
over all set-valued tableaux $S$ of shape $\lambda/\mu$.  The proof is
based on a bijection between set-valued tableaux and monomials
equivalent to compatible pairs.

Given an integer $n$ such that $\pi_{\lambda/\mu} \in S_n$, we can
formulate this as a bijection between set-valued tableaux of shape
$\lambda/\mu$ and (diagrams of) compatible pairs $(a,i)$ with $w(a) =
\pi_0 \pi_{\lambda/\mu}^{-1} \pi_0$, where $\pi_0$ is the longest
permutation in $S_n$.  Number the north-west to south-east diagonals
of $\lambda$ (and $\mu$) consecutively from right to left, so that the
upper-left box of $\lambda$ is in diagonal number $n-k$.  Then the
set-valued tableau $S$ of shape $\lambda/\mu$ is mapped to the diagram
in which column $p$ consists of the diagonal numbers of the boxes of
$S$ that contain the integer $p$.

\begin{example} \label{E:setval_diag}
If $n-k=4$, then the bijection makes the assignments
\[ \tableau{&&{1}&{23}\\&{12}&{24}\\{2}&{36}&{7}}
   \ \mapsto \ 
   \raisebox{13.5pt}{\begin{picture}(0,0)
     \put(11,11){\line(1,0){55}}
     \end{picture}
     \tabl{{2}&{1}&{1}&{3}&&{5}&{4}\\{4}&{3}&{5}\\&{4}\\&{6}}}
   \ \ \ \ \text{and} \ \ \ \ 
   \tableau{&&{1}&{1}\\&{12}&{2}\\{12}&{3}&{3}}
   \ \mapsto \ 
   \tableau{{1}&{3}&{4}\\{2}&{4}&{5}\\{4}&{6}\\{6}}
   \,.
\]
\end{example}

Recall that the column reading word of a set-valued tableau is
obtained by reading its boxes from bottom to top and then left to
right.  The integers of each box are arranged in increasing order.
The set-valued tableaux in Example~\ref{E:setval_diag} have column
reading words $2\,3\,6\,1\,2\,7\,2\,4\,1\,2\,3$ and
$1\,2\,3\,1\,2\,3\,2\,1\,1$.  Recall also that a word is a {\em
  reverse lattice word\/} if every occurrence of an integer $i$ with
$i > 1$ is followed by more $i-1$'s than $i$'s.  The {\em content\/}
of a word is the integer sequence $(\nu_1,\nu_2,\dots)$ where $\nu_i$
is the number of occurrences of $i$ in the word.  The following lemma
says that the condition that the diagram of a compatible pair is an
increasing tableau naturally generalizes the condition that the
reading word of a tableau is a reverse lattice word.

\begin{lemma} \label{L:lr_incr}
  If the set-valued tableau $S$ is mapped to the diagram $T$, then the
  column reading word of $S$ is a reverse lattice word if and only if
  $T$ is an increasing tableau.
\end{lemma}
\begin{proof}
  Consider an integer $i > 1$ contained in some box $B$ of $S$.  The
  lemma follows from the observation that all occurrences of $i$ and
  $i-1$ that follow the integers of $B$ in the column reading word of
  $S$ have diagonal numbers that are strictly smaller than the
  diagonal number of $B$, and all other occurrences of $i$ and $i-1$
  have diagonal numbers that are larger than or equal to the diagonal
  number of $B$.
\end{proof}

\begin{cor}[Theorem 6.9 of \cite{buch:littlewood-richardson}]
  The coefficient $c_{\pi_{\lambda/\mu},\nu}$ is equal to the number
  of set-valued tableaux $S$ of shape $\lambda/\mu$ such that the
  column reading word of $S$ is a reverse lattice word with content
  $\nu$.
\end{cor}
\begin{proof}
  Notice that if the set-valued tableau $S$ is mapped to the diagram
  $T$, then the content of $S$ is the list of column lengths of $T$.
  The number of set-valued tableau $S$ of the corollary is therefore
  equal to the number of increasing tableaux $T$ of shape conjugate to
  $\nu$ such that $w(T) = \pi_0 \pi_{\lambda/\mu}^{-1} \pi_0$, which
  by Theorem~\ref{T:stabcoef} equals $c_{\pi_0 \pi_{\lambda/\mu}
    \pi_0, \nu'} = c_{\pi_{\lambda/\mu}, \nu}$ (see
  \cite[Lemma~3.4]{buch:littlewood-richardson}).
\end{proof}

\subsection{Generating increasing tableaux}

For practical applications of Theorem~\ref{T:stabcoef}, it is
desirable to have an efficient algorithm for generating all increasing
tableaux which represent a given permutation.  We will address the
following slightly more general problem.

\begin{problem}
  Given a column $C_0$ of boxes containing integers that increase from
  top to bottom and a permutation $\pi \in S_\infty$, find all
  increasing tableaux $T$ such that $w(T) = \pi$ and so that $T$ can
  be attached to the right hand side of $C_0$ to form a larger
  increasing tableau.
\end{problem}

We can generate the solutions to this problem as follows.  First we
find all pairs $(C, \sigma)$ consisting of an increasing column $C$
and a permutation $\sigma$, such that $\pi = w(C) \cdot \sigma$ and
$C$ can be attached to the right hand side of $C_0$.  For each such
pair we recursively find all increasing tableaux $T'$ for which $w(T')
= \sigma$ and $T'$ can be attached to the right side of $C$, thus
forming one of the solutions $T$.  Notice that it is sufficient to
consider pairs $(C,\sigma)$ such that $\sigma(i)=i$ for the smallest
integer $i$ for which $\pi(i) > i$.  Furthermore, these pairs can be
generated very quickly.

However, the algorithm in its present form is not efficient, because
in many cases there are no increasing tableaux $T$ satisfying the
stated conditions, and it may require many recursive applications of
the algorithm to discover this.  We will fix this problem by
describing an easy way to decide up front if at least one solution $T$
exists.

Stanley has proved \cite{stanley:on} that the Schur expansion of his
symmetric function $F_{\pi^{-1}}$ contains two special terms, each
with coefficient one, which are indexed by partitions that are minimal
and maximal in the dominance order among partitions occurring in
$F_{\pi^{-1}}$.  Accordingly, there is exactly one increasing tableau
representing $\pi$ on each of these shapes.  Let $M_\pi$ be the unique
increasing tableau of the maximal shape.  Then the integers $i_1 < i_2
< \dots < i_p$ in the top row of $M_\pi$ satisfy that each $i_k$ is
the largest descent position smaller than $i_{k+1}$ of the permutation
$\pi s_{i_p} s_{i_{p-1}} \cdots s_{i_{k+1}}$ (we set
$i_{p+1}=\infty$), and the permutation $\pi s_{i_p} \cdots s_{i_1}$
has no descent positions smaller than $i_1$.  This characterizes the
top row of $M_\pi$, and the part of $M_\pi$ below this row is equal to
$M_{\pi s_{i_p} \cdots s_{i_1}}$.

We leave it as an exercise for the reader to show that the integers of
the first column of $M_\pi$ are larger than the integers in the first
column of any other increasing tableau representing $\pi$.  In other
words, if $w(T) = \pi$ and the leftmost box of row $r$ in $T$ contains
the integer $x$, then either $M_\pi$ has fewer than $r$ rows, or the
first integer in its $r$th row is larger than or equal to $x$.  This
property of $M_\pi$ implies that the set of solutions $T$ to our
problem is nonempty if and only if $M_\pi$ is a solution.  When this
criterion is incorporated, our algorithm is fairly efficient.

\subsection{Products of increasing tableaux} \label{SS:product}

Given two increasing tableau $T_1$ and $T_2$, we let $T_1 \cdot T_2$
denote the increasing tableau obtained by Hecke inserting the word of
$T_1$ into $T_2$.  More precisely, if $a_1 a_2 \cdots a_p$ is the word
of $T_1$ then we define $T_1 \cdot T_2 = (a_1 \hcins (a_2\hcins
(\dotsm (a_p \hcins T_2)\dotsm)))$.  When the concatenation of the
words of $T_1$ and $T_2$ is a reduced word of a permutation, this
product agrees with the Coxeter-Knuth product, which is known to be
associative.  Unfortunately our more general product of increasing
tableaux is not associative.  We make the convention that $T_1\cdot
T_2\cdot T_3$ means $T_1\cdot (T_2\cdot T_3)$.

\begin{example} Let $T_1=\tableau{{1}}\,$, $T_2=\tableau{{1}&{5}\\{4}}\,$, 
and $T_3=\tableau{{2}}\,$.  Then
\begin{equation*}
(T_1\cdot T_2) \cdot T_3 \ = \ \tableau{{1}&{4}&{5}\\{4}}\,\cdot
\tableau{{2}} \ = \ \tableau{{1}&{2}&{5}\\{4}}
\end{equation*}
whereas
\begin{equation*}
T_1 \cdot (T_2\cdot T_3) \ = \ \tableau{{1}} \cdot \tableau{{1}&{2}\\{4}&{5}} 
\ = \ \tableau{{1}&{2}&{5}\\{4}&{5}} \ .
\end{equation*}
\end{example}

The product of increasing tableaux has the following properties, whose
proofs are straightforward from the definitions.

\begin{lemma} \label{L:tabprod}
Let $T$ and $T'$ be increasing tableaux.  Then we have
\begin{enumerate}
\item $w(T \cdot T') = w(T) \cdot w(T')$.
\item Suppose $T$ is cut along a vertical line into $T_\lft$ and
  $T_\rgt$.  Then $T = T_\lft \cdot T_\rgt$.
\item Suppose $T$ is cut along a horizontal line into tableaux
  $T_\bott$ and $T_\topp$.  Then $T = T_\bott \cdot T_\topp$.
\end{enumerate}
\end{lemma}

\subsection{Decreasing tableaux} \label{SS:decrtab}

For our applications to $K$-theoretic factor sequences in the next
section, it is more natural to work with decreasing tableaux, which by
definition are Young tableaux with strictly decreasing rows and
columns.  By regarding decreasing tableaux as increasing tableaux with
the order of the natural numbers inverted, we obtain well defined
operations of Hecke insertion and products of decreasing tableaux as
in Sections \ref{SS:heckeins} and \ref{SS:product}.  For example we
have
\[ \tableau{{5}&{3}\\{4}} \ \cdot \ \tableau{{6}&{3}&{1}\\{4}&{2}\\{3}\\{2}} 
   \ = \ \tableau{{6}&{4}&{3}&{1}\\{5}&{3}&{2}\\{4}\\{2}} 
\]

If $T$ is a decreasing tableau, we let $w(T)$ be the unique
permutation that has the column word of $T$ as a Hecke word.
Lemma~\ref{L:tabprod} then remains true for products of decreasing
tableaux.  In addition we have the following decreasing version of
Theorem~\ref{T:stabcoef}.

\begin{thmp} \label{T:stabcoefd}
  For any permutation $\pi$, the coefficient $c_{\pi,\lambda}$ of
  (\ref{E:stabcoef}) equals $(-1)^{|\lambda|-\ell(\pi)}$ times the
  number of decreasing tableaux $T$ of shape $\lambda$ such that
  $w(T)=\pi$.
\end{thmp}
\begin{proof}
  Let $\pi \in S_n$ and let $\pi_0 \in S_n$ be the longest permutation.
  By replacing each entry $x$ in a tableau with $n-x$, we obtain a
  bijection between the decreasing tableaux representing $\pi$ and the
  increasing tableaux representing $\pi_0 \pi \pi_0$.  The theorem
  therefore follows from the identity $G_{\pi_0 \pi^{-1} \pi_0} = G_\pi$,
  which is a consequence of Fomin and Kirillov's construction of
  Grothendieck polynomials.
\end{proof}

\section{Quiver coefficients}
\label{S:quiver}

\subsection{$K$-theoretic factor sequences}

Let $r = \{ r_{ij} \}$ be a set of rank conditions for $0 \leq i \leq
j \leq n$, and set $N = e_0 + \dots + e_n$ where $e_i = r_{ii}$.  A
result of Zelevinsky shows that when the base variety ${\mathfrak X}$
is a product of matrix spaces, the quiver variety $\Omega_r \subset
{\mathfrak X}$ is isomorphic to a dense open subset of a Schubert
variety \cite{zelevinsky:two}. The {\em Zelevinsky permutation\/}
corresponding to this Schubert variety was used in
\cite{knutson.miller.ea:four} to prove the ratio formula for quiver
varieties.

With the notation from \cite{buch:alternating}, the Zelevinsky
permutation can be constructed as a product of permutations as follows
(see \cite[Prop.~1.6]{knutson.miller.ea:four} for a different
construction).  Extend the rank conditions $r = \{ r_{ij} \}$ by
setting $r_{ij} = e_j + \dots + e_i$ for $0 \leq j < i \leq n$.  Then
define decreasing tableaux $U_{ij}$ as in the introduction, but for
all $0 \leq i < n$ and $0 < j \leq n$.  The corresponding permutations
$W_{ij} = w(U_{ij})$ are given by
\[ W_{ij}(p) = \begin{cases}
  p + r_{i,j-1} - r_{ij} & \text{if $r_{ij} < p \leq r_{i+1,j}$} \\
  p - r_{i+1,j} + r_{ij} & \text{if $r_{i+1,j} < p \leq
    r_{i+1,j} + r_{i,j-1} - r_{ij}$} \\
  p & \text{otherwise.}
\end{cases}\]
The Zelevinsky permutation can now be defined by $z(r) = \prod_{j=1}^n
\prod_{i=0}^{n-1} W_{ij}$.  The descent positions of $z(r)$ are
contained in the set $\{ r_{nj} : 0 < j \leq n\}$, and the descent
positions of $z(r)^{-1}$ are contained in $\{ r_{i0} : 0 \leq i < n
\}$.

For each $1 \leq j \leq n-1$ we set $\delta_j = W_{jj} W_{j+1,j}
\cdots W_{n-1,j} \in S_N$.  A {\em KMS-factorization\/} for the rank
conditions $r$ is any sequence $(\pi_1,\dots,\pi_n)$ of permutations
with $\pi_i \in S_{e_{i-1}+e_i}$, such that the Zelevinsky permutation
$z(r)$ is equal to the Hecke product
\[ \pi_1 \cdot \delta_1 \cdot \pi_2 \cdot \delta_2 \cdots \delta_{n-1}
   \cdot \pi_n \,.
\]
These sequences of permutations generalize the notion of a minimal
lace diagram from \cite{knutson.miller.ea:four} and give the index set
in the $K$-theoretic stable component formula \eqref{E:stabcomp} from
\cite{buch:alternating, miller:alternating}.

We begin by defining a {\em $K$-theoretic factor sequence\/} for the
rank conditions $r$ to be any sequence $(T_1,\dots,T_n)$ of decreasing
tableaux, such that $(w(T_1), \dots, w(T_n))$ is a KMS-factorization
for $r$.  With this definition, Theorem~\ref{T:quivcoef} is an
immediate consequence of Theorem~\ref{T:stabcoef} combined with the
$K$-theoretic stable component formula \eqref{E:stabcomp}.  To obtain
the inductive definition of factor sequences given before Theorem
\ref{T:quivcoef} we need the following result, proved in
\cite[Thm.~7]{buch:alternating}, which shows that KMS-factorizations
can themselves be defined as `factor sequences'.  Recall the definition
of $\wb r$ from Section~\ref{S:intro_quiver}.

\begin{thm} \label{T:kmsfacseq} (a) If $(\pi_1,\dots,\pi_n)$ is a
  KMS-factorization for $r$, then each permutation $\pi_i$ has a
  reduced factorization $\pi_i = \rho_{i-1} \cdot W_{i-1,i} \cdot
  \sigma_i$ with $\rho_{i-1} \in S_{e_{i-1}}$ and $\sigma_i \in
  S_{e_i}$, such that $\rho_0 = \sigma_n = 1$.
  
  (b) Let $\sigma_1,\rho_1,\dots,\sigma_{n-1},\rho_{n-1}$ be
  permutations with $\sigma_i,\rho_i \in S_{e_i}$.  Then the sequence
  $(W_{01}\cdot \sigma_1, \rho_1 \cdot W_{12} \cdot \sigma_2, \dots,
  \rho_{n-1} \cdot W_{n-1,n})$ is a KMS-factorization for $r$ if and
  only if $(\sigma_1 \cdot \rho_1, \sigma_2 \cdot \rho_2, \dots,
  \sigma_{n-1} \cdot \rho_{n-1})$ is a KMS-factorization for $\wb r$.
\end{thm}

We also need the following statement.

\begin{lemma} \label{L:upperleft} Let $T$ be any decreasing tableau
  such that $w(T) \in S_m$, and for some integers $a,b < m$ we have
  $w(T)(p) \leq b$ for all $a < p \leq m$.  Then $T$ contains the
  rectangle $R = (m-b)\times (m-a)$ in its upper left corner.  The
  upper-left box of $R$ equals $m-1$, and the boxes of $R$ decrease by
  one for each step down or to the right.
\end{lemma}
\begin{proof}
  After deleting the contents of some boxes of $T$, the permutation
  $w(T)$ becomes equal to the south-west to north-east (permutation)
  product of the simple transpositions corresponding to the non-empty
  boxes in $T$.  Since the integers in all these boxes are smaller
  than $m$, the assumption that $w(T)(m) \leq b$ implies that the top
  of the first column of $T$ must contain the integers $m-1, m-2,
  \dots, b$.  The assumption that $w(T)(m-1) \leq b$ then implies that
  the second column of $T$ starts with $m-2,m-3,\dots,b-1$, etc.
\end{proof}

Let $(T_1,T_2) \mapsto T_1 \cdot T_2$ be the product of decreasing
tableaux from Section \ref{SS:decrtab}.

\begin{cor}
  A sequence of decreasing tableaux $(T_1,\dots,T_n)$ is a
  $K$-theoretic factor sequence for the rank conditions $r$ if and
  only if there exist decreasing tableaux $A_i, B_i$ for $1 \leq i
  \leq n-1$, such that $T_i = B_{i-1} \cdot U_{i-1,i} \cdot A_i$ for
  each $i$ (with $B_0 = A_n = \emptyset$) and $(A_1 \cdot B_1, \dots,
  A_{n-1}\cdot B_{n-1})$ is a $K$-theoretic factor sequence for $\wb
  r$.
\end{cor}
\begin{proof}
  Let $(T_1,\dots,T_n)$ be a factor sequence for $r$ and
  $(\pi_1,\dots,\pi_n)$ the corresponding KMS-factorization.  It
  follows from Theorem~\ref{T:kmsfacseq} (a) that $\pi_i \in
  S_{e_{i-1} + e_i - r_{i-1,i}}$ and that $\pi_i(p) \leq e_{i-1}$ for
  all $e_i < p \leq e_{i-1}+e_i-r_{i-1,i}$.  Since $T_i$ represents
  $\pi_i$, Lemma~\ref{L:upperleft} implies that $T_i$ contains the
  tableau $U_{i-1,i}$ in its upper-left corner.  Now write $T_i =
  B_{i-1} \cdot U_{i-1,i} \cdot A_i$ where $A_i$ is the part of $T_i$
  to the right of $U_{i-1,i}$ and $B_{i-1}$ is the part below
  $U_{i-1,i}$ and $A_i$.
  \[ T_i \ \ = \ \ \raisebox{-25pt}{\pic{50}{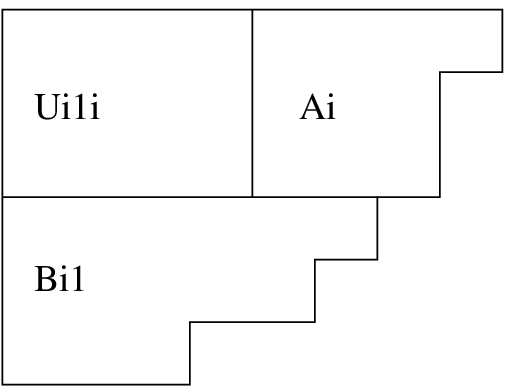}} \]
  Then we have $\pi_i = w(B_{i-1}) \cdot W_{i-1,i} \cdot w(A_i)$ by
  Lemma~\ref{L:tabprod}, and all entries of $A_i$ and $B_i$ are
  smaller than $e_i$.  Since all descent positions of $z(r) = \pi_1
  \cdot \delta_1 \cdot \pi_2 \cdots \delta_{n-1} \cdot \pi_n$ are
  greater than or equal to $e_n$, the same must be true for $\pi_n$,
  so $A_n$ must be empty.  Similarly, since the descent positions of
  $\pi_1^{-1}$ are greater than or equal to $e_0$, $B_0$ is empty.
  Now it follows from Theorem~\ref{T:kmsfacseq} (b) that $(w(A_1)\cdot
  w(B_1), \dots, w(A_{n-1})\cdot w(B_{n-1}))$ is a KMS-factorization
  for $\wb r$, or equivalently that $(A_1 \cdot B_1, \dots, A_{n-1}
  \cdot B_{n-1})$ is a factor sequence.

  On the other hand, if we are given decreasing tableaux $A_1,B_1,
  \dots, A_{n-1},B_{n-1}$ such that $(A_1\cdot B_1, \dots,
  A_{n-1}\cdot B_{n-1})$ is a factor sequence for $\wb r$ then
  Theorem~\ref{T:kmsfacseq} (a) implies that the entries of $A_i$ and
  $B_i$ are smaller than $r_{i-1,i}+r_{i,i+1}-r_{i-1,i+1} \leq e_i$,
  so it follows from Theorem \ref{T:kmsfacseq} (b) that $(U_{01} \cdot
  A_1, B_1 \cdot U_{12} \cdot A_2, \dots, B_{n-1} \cdot U_{n-1,n})$ is
  a factor sequence for $r$.
\end{proof}

This completes the proof of Theorems \ref{T:quivcoef} and \ref{T:fackms}.

\begin{example} \label{X:facseq}
  Consider a sequence of vector bundles $E_0 \to E_1 \to E_2 \to E_3$
  of ranks $1, 4, 3, 3$ together with the rank conditions
\[ r \ = \ 
   \left\{ \begin{matrix} 
     r_{00} && r_{11} && r_{22} && r_{33} \\
     & r_{01} && r_{12} && r_{23} \\
     && r_{02} && r_{13} \\
     &&& r_{03}
  \end{matrix} \right\} 
  \ = \ 
   \left\{ \begin{matrix} 
     1 && 4 && 3 && 3 \\
     & 1 && 2 && 2 \\
     && 1 && 1 \\
     &&& 0
  \end{matrix} \right\} \,.
\]
These rank conditions result in the following diagram of decreasing
tableaux $U_{ij}$:
\[ \begin{matrix}
  \,\emptyset\, && \tableau{{4}&{3}} && \tableau{{3}} \vspace{1mm}\\
  & \,\emptyset\, && \tableau{{2}} \vspace{1mm}\\
  && \tableau{{1}}
  \end{matrix}
\]
The two bottom rows of this diagram produce the following three factor
sequences for the inductive rank conditions $\wb r$:
\[ \left(\,\tableau{{1}} \,,\, \tableau{{2}}\, \right) \,,\, 
   \left(\emptyset \,,\, \tableau{{2}\\{1}}\, \right) \,,\, 
   \left(\,\tableau{{1}} \,,\, \tableau{{2}\\{1}}\, \right) \,. 
\]
Since the decreasing tableau $\tableau{{2}\\{1}}$ has the
factorizations
\[ \tableau{{2}\\{1}} 
   \ = \ \tableau{{2}\\{1}} \cdot \emptyset
   \ = \ \tableau{{1}} \cdot \tableau{{2}}
   \ = \ \emptyset \cdot \tableau{{2}\\{1}}
   \ = \ \tableau{{2}\\{1}} \cdot \tableau{{2}}
   \ = \ \tableau{{1}} \cdot \tableau{{2}\\{1}}
\]
and a single box has the factorizations $\tableau{{x}} = \tableau{{x}}
\cdot \emptyset = \emptyset \cdot \tableau{{x}} = \tableau{{x}} \cdot
\tableau{{x}}$\,, the factor sequences for $\wb r$ produce the following
list of factor sequences for $r$:
\vspace{2mm}

\noindent$\mbox{}\hspace{13mm}\,
\left(\,\tableau{{1}}\,,\,\tableau{{4}&{3}&{2}}\,,\,\tableau{{3}}\,\right),
\left(\tableau{{1}}\,,\,\tableau{{4}&{3}}\,,\,\tableau{{3}\\{2}}\,\right),
\left(\emptyset,\,\tableau{{4}&{3}&{2}\\{1}}\,,\,\tableau{{3}}\,\right),
\vspace{1mm}\\\mbox{}\hspace{13mm}
\left(\emptyset,\,\tableau{{4}&{3}\\{1}}\,,\,\tableau{{3}\\{2}}\,\right),
\left(\emptyset,\,\tableau{{4}&{3}&{1}}\,,\,\tableau{{3}\\{2}}\,\right),
\left(\emptyset,\,\tableau{{4}&{3}}\,,\,\tableau{{3}\\{2}\\{1}}\,\right),
\vspace{1mm}\\\mbox{}\hspace{13mm}
\left(\tableau{{1}}\,,\,\tableau{{4}&{3}&{2}\\{1}}\,,\,\tableau{{3}}\,\right),
\left(\tableau{{1}}\,,\,\tableau{{4}&{3}\\{1}}\,,\,\tableau{{3}\\{2}}\,\right),
\left(\tableau{{1}}\,,\,\tableau{{4}&{3}&{2}}\,,\,\tableau{{3}\\{2}}\,\right),
\vspace{1mm}\\\mbox{}\hspace{13mm}
\left(\tableau{{1}}\,,\,\tableau{{4}&{3}&{1}}\,,\,\tableau{{3}\\{2}}\,\right),
\left(\emptyset,\,\tableau{{4}&{3}&{2}\\{1}}\,,\,\tableau{{3}\\{2}}\,\right),
\left(\emptyset,\,\tableau{{4}&{3}&{1}\\{1}}\,,\,\tableau{{3}\\{2}}\,\right),
\vspace{1mm}\\\mbox{}\hspace{13mm}
\left(\emptyset,\,\tableau{{4}&{3}&{1}}\,,\,\tableau{{3}\\{2}\\{1}}\,\right),
\left(\tableau{{1}}\,,\,\tableau{{4}&{3}}\,,\,\tableau{{3}\\{2}\\{1}}\,\right),
\left(\emptyset,\,\tableau{{4}&{3}\\{1}}\,,\,\tableau{{3}\\{2}\\{1}}\,\right),
\vspace{1mm}\\\mbox{}\hspace{13mm}
\left(\tableau{{1}}\,,\,\tableau{{4}&{3}&{2}\\{1}}\,,\,\tableau{{3}\\{2}}\,\right),
\left(\tableau{{1}}\,,\,\tableau{{4}&{3}&{1}\\{1}}\,,\,\tableau{{3}\\{2}}\,\right),
\left(\tableau{{1}}\,,\,\tableau{{4}&{3}&{1}}\,,\,\tableau{{3}\\{2}\\{1}}\,\right),
\vspace{1mm}\\\mbox{}\hspace{13mm}
\left(\emptyset,\,\tableau{{4}&{3}&{1}\\{1}}\,,\,\tableau{{3}\\{2}\\{1}}\,\right),
\left(\tableau{{1}}\,,\,\tableau{{4}&{3}\\{1}}\,,\,\tableau{{3}\\{2}\\{1}}\,\right),
\left(\tableau{{1}}\,,\,\tableau{{4}&{3}&{1}\\{1}}\,,\,\tableau{{3}\\{2}\\{1}}\,\right)
$
\vspace{2mm}\\
These factor sequences can also be obtained by first working out the
13 possible KMS-factorizations for $r$, for example by using Theorem
\ref{T:kmsfacseq} or the transformations on KMS-factorizations given
in \cite[\S 5]{buch.feher.ea:positivity}.  
We conclude that the Grothendieck class of the quiver variety
$\Omega_r(E_\bull)$ is obtained by replacing each tensor $G_{\mu_1}
\otimes G_{\mu_2} \otimes G_{\mu_3}$ in the following expression with
the class $G_{\mu_1}(E_1-E_0) \cdot G_{\mu_2}(E_2-E_1) \cdot
G_{\mu_3}(E_3-E_2)$:
\[\begin{split}
&     G_{1} \otimes G_{3} \otimes G_{1} 
   + G_{1} \otimes G_{2} \otimes G_{11} 
   + 1 \otimes G_{31} \otimes G_{1} 
   + 1 \otimes G_{21} \otimes G_{11} 
\\&
   + 1 \otimes G_{3} \otimes G_{11} 
   + 1 \otimes G_{2} \otimes G_{111} 
   - G_{1} \otimes G_{31} \otimes G_{1}
   - G_{1} \otimes G_{21} \otimes G_{11}
\\&
   - 2 \cdot G_{1} \otimes G_{3} \otimes G_{11}
   - 2 \cdot 1 \otimes G_{31} \otimes G_{11}
   - 1 \otimes G_{3} \otimes G_{111}
   - G_{1} \otimes G_{2} \otimes G_{111} 
\\&
   - 1 \otimes G_{21} \otimes G_{111}
   + 2 \cdot G_{1} \otimes G_{31} \otimes G_{11}
   + G_{1} \otimes G_{3} \otimes G_{111}
   + 1 \otimes G_{31} \otimes G_{111}
\\&
   + G_{1} \otimes G_{21} \otimes G_{111}
   - G_{1} \otimes G_{31} \otimes G_{111} \,.
\end{split}\]
\end{example}

\subsection{Grothendieck polynomials for Zelevinsky permutations}

Using results about Demazure characters it was proved in
\cite{knutson.miller.ea:four} that cohomological quiver coefficients
are special cases of the Stanley coefficients associated to the
Zelevinsky permutation $z(r)$.  As an application of our results, we
will prove more generally that the $K$-theoretic quiver coefficients
are special cases of the coefficients $c_{z(r),\lambda}$ in the
expansion \eqref{E:stabcoef} of the stable Grothendieck polynomial for
$z(r)$.  This result also sharpens the fact from
\cite{buch:alternating, buch.sottile.ea:quiver} that quiver
coefficients are special cases of the decomposition coefficients of
Grothendieck polynomials studied in \cite{buch.kresch.ea:grothendieck}
(see Section \ref{SS:decomp}).  Given a sequence of partitions $\mu =
(\mu_1,\dots,\mu_n)$ such that $\mu_i$ is contained in the rectangle
$e_i \times e_{i-1}$, let $\lambda(\mu)$ be the partition obtained by
concatenating the partitions $(e_0+e_1+\dots+e_{i-2})^{e_i} + \mu_i$
for $i = n, n-1, \dots, 1$.

\begin{thm} \label{T:quivstab}
  For any set of rank conditions $r$ and sequence of partitions $\mu$
  we have $c_\mu(r) = c_{z(r),\lambda(\mu)}$.
\end{thm}

Our proof of the above identity is based on a bijection between the
$K$-theoretic factor sequences for $r$ and the decreasing tableaux
representing $z(r)$.  Given a sequence $(T_1,\dots,T_n)$ of decreasing
tableaux, such that each tableau $T_i$ can be contained in the
rectangle $e_i \times e_{i-1}$ and all entries of $T_i$ are smaller
than $e_{i-1}+e_i$, we let $\Phi(T_1,\dots,T_n)$ denote the decreasing
tableau constructed from this sequence as well as the tableaux
$U_{ij}$ for $i \geq j$ as follows.
\[ \Phi(T_1,\dots,T_n) \ \ = \ \ \raisebox{-52pt}{\pic{60}{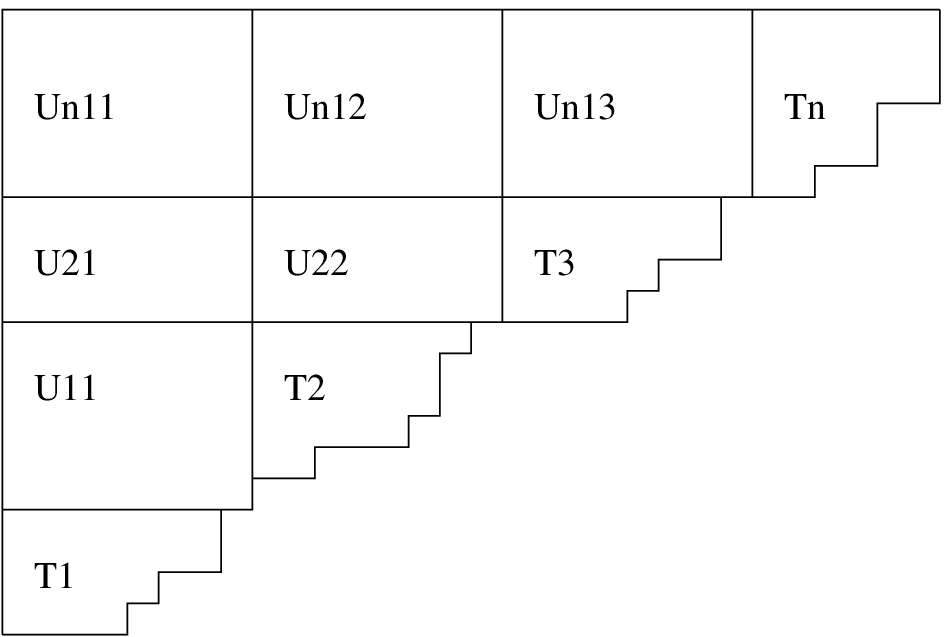}} \]
Notice that the upper-left box of $U_{n-1,1}$ is equal to $N-1$, and
the boxes in the union of tableaux $U_{ij}$ decrease by one for each
step down or to the right.  Theorem~\ref{T:quivstab} follows from the
following proposition combined with Theorems~\ref{T:stabcoef} and
\ref{T:quivcoef}.

\begin{prop}
  The map $(T_1,\dots,T_n) \mapsto \Phi(T_1,\dots,T_n)$ gives a
  bijection of the set of all $K$-theoretic factor sequences for $r$
  with the set of all decreasing tableaux representing $z(r)$.
\end{prop}
\begin{proof}
  Since the permutation of a decreasing tableau can be defined as the
  south-west to north-east Hecke product of the simple reflections
  given by the boxes of the tableau, it follows from the definition of
  KMS-factorizations that $(T_1,\dots,T_n)$ is a factor sequence if
  and only if $\Phi(T_1,\dots,T_n)$ represents the Zelevinsky
  permutation $z(r)$.  It remains to show that any decreasing tableau
  $T$ representing $z(r)$ contains the arrangement of rectangular
  tableaux $U_{ij}$ in its upper-left corner, and has no boxes
  strictly south-east of the tableaux $U_{ii}$ for $1 \leq i \leq
  n-1$.  The inclusion of the tableaux $U_{ij}$ in $T$ follows from
  Lemma~\ref{L:upperleft} because $z(r) \in S_N$ and for each $0 < i <
  n$ and $p > r_{ni}$ we have $z(r)(p) \leq r_{i0}$, see
  \cite[Prop~1.6]{knutson.miller.ea:four} or \cite[Lemma
  3.1]{buch:alternating}.

  To see that $T$ contains no boxes strictly south-east of $U_{ii}$,
  we use that the \linebreak Grothendieck polynomial $\Groth_{\wh
    z(r)}(x_1,\dots,x_N)$ is separately symmetric in each group of
  variables $\{ x_p \mid r_{n,i} < p \leq r_{n,i-1} \}$, where $\wh
  z(r) = \pi_0^{(N)} z(r)^{-1} \pi_0^{(N)}$ and $\pi_0^{(N)}$ is the longest
  permutation in $S_N$.  This is true because the descent positions of
  $\wh z(r)$ are contained in the set $\{r_{nj}\mid 0 < j \leq n\}$.
  It follows that the exponent of $x_{r_{ni}+1}$ in any monomial of
  $\Groth_{\wh z(r)}(x_1,\dots,x_N)$ is less than or equal to $N -
  r_{n,i-1} = r_{i-2,0}$.  Now $T$ can be used to construct a unique
  compatible pair $(a,k)$ for $\wh z(r)$, such that $T$ contains the
  integer $p$ in some box of row $q$ if and only if $(a_l,k_l) =
  (N-p,q)$ for some $l$.  Since this pair contributes the monomial
  $x^k$ to $\Groth_{\wh z(r)}(x_1,\dots,x_N)$, it follows that row
  $r_{ni}+1$ of $T$ has at most $r_{i-2,0}$ boxes.  This means exactly
  that $T$ contains no boxes south-east of $U_{i-1,i-1}$, as required.
\end{proof}


\setcounter{MaxMatrixCols}{11}
\begin{example}
  The Zelevinsky permutation for the rank conditions $r$ of
  Example~\ref{X:facseq} is given by $z(r) = \begin{pmatrix}
    1&2&3&4&5&6&7&8&9&10&11 \\ 2&6&9&1&10&11&3&4&7&8&5 \end{pmatrix}
  \in S_{11}$.  The decreasing tableaux representing this permutation
  are obtained by attaching the factor sequences for $r$ to the bottom
  side, the middle corner, and the right side of the tableau
\[ \tableau{{10}&{9}&{8}&{7}&{6}\\
{9}&{8}&{7}&{6}&{5}\\
{8}&{7}&{6}&{5}&{4}\\
{7}\\
{6}\\
{5}}
 \ . \]
\end{example}

\subsection{Universal Grothendieck polynomials}

\label{SS:decomp} Fulton's \textit{universal Schubert polynomials}
\cite{fulton:universal} describe certain quiver varieties associated
to a sequence of vector bundles $E_1 \to \dots \to E_{n-1} \to E_n \to
F_n \to F_{n-1} \to \dots \to F_1$ over $X$, such that $\rank(E_i) =
\rank(F_i) = i$ for each $i$.  They are also known to specialize to,
e.g., the quantum Schubert polynomials
\cite{fomin.gelfand.postnikov:quantum}, where a nonrecursive
combinatorial formula was given in \cite{buch.kresch.ea:schubert}. We
now describe an extension of this result to $K$-theory.

Given a permutation $\pi \in S_{n+1}$, we let $\Omega_\pi \subset X$ be
the degeneracy locus of points where the rank of each composed map
$E_q \to F_p$ is at most equal to the number of integers $i \leq p$ 
such that $\pi(i) \leq q$.  The quiver formula \eqref{E:kquiver}
can be applied to give a formula
\begin{equation} \label{E:univgroth}
  [\OO_{\Omega_\pi}] = \sum_\mu c^{(n)}_{\pi,\mu}\,
  G_{\mu_1}(E_2-E_1) \cdots G_{\mu_n}(F_n-E_n) \cdots
  G_{\mu_{2n-1}}(F_1-F_2)
\end{equation}
for the Grothendieck class of $\Omega_\pi$, where the coefficients
$c^{(n)}_{\pi,\mu}$ are special cases of quiver coefficients.  It was
shown in \cite{buch:grothendieck} that the coefficients
$c_{\pi,\lambda}$ of the expansion \eqref{E:stabcoef} of the stable
Grothendieck polynomial for $\pi$ can be obtained as the
specializations
$c^{(n)}_{\pi,(\emptyset^{n-1},\lambda,\emptyset^{n-1})}$, where
$\emptyset^{n-1}$ denotes a sequence of $n-1$ empty partitions. More
generally, it was proved in \cite[Thm.~4]{buch.kresch.ea:grothendieck}
that the coefficients $c^{(n)}_{\pi,\lambda}$ can be used to expand a
double Grothendieck polynomial as a linear combination of products of
stable Grothendieck polynomials applied to disjoint intervals of
variables.  In \cite{buch.kresch.ea:grothendieck}, the formula
\eqref{E:univgroth} was also used to prove that
\[ [\OO_{\Omega_\pi}] = \sum (-1)^{\ell(\sigma_1\cdots \sigma_{2n-1} \pi)}
   G_{\sigma_1}(E_2-E_1) \cdots G_{\sigma_n}(F_n-E_n) \cdots
   G_{\sigma_{2n-1}}(F_1-F_2)
\]
where this sum is over all sequences of permutations
$(\sigma_1,\dots,\sigma_{2n-1})$ such that $\sigma_i \in
S_{\min(i,2n-i)+1}$ and $\pi$ is equal to the Hecke product $\sigma_1
\cdot \sigma_2 \cdots \sigma_{2n-1}$.  Combining this with
Theorem~\ref{T:stabcoef}, we obtain the following generalization of
\cite[Thm.~1]{buch.kresch.ea:schubert}.

\begin{thm}
\label{T:univgroth}
The coefficient $c^{(n)}_{\pi,\mu}$ of \eqref{E:univgroth} is equal to
$(-1)^{\sum |\mu_i| - \ell(\pi)}$ times the number of sequences
$(T_1,\dots,T_{2n-1})$ of increasing tableaux of shapes $(\mu_1,
\dots, \mu_{2n-1})$, such that the entries of $T_i$ are at most
$\min(i,2n-i)$ and $w(T_{2n-1} \dotsm T_2 \cdot T_1)=\pi^{-1}$.
\end{thm}

\section*{Acknowledgments}

We would like to thank S.~Fomin, M.~Haiman, A.~Postnikov and J.~Remmel
for helpful conversations.  This work was partially completed while
AB, MS and AY were in residence together at the Park City Mathematics
Institute, during the program on ``Geometric combinatorics'' during
July 2004.  AK was partially supported by an EPSRC Advanced Research
Fellowship.  MS was supported in part by NSF grant DMS-0401012.  HT
was supported in part by NSF grant DMS-0401082.  AY thanks NSERC for
providing support as a visitor to the Fields Institute in Toronto,
during the 2005 semester on ``The Geometry of String Theory''; he
would also like to thank the NSF for supporting a visit to the
Mittag-Leffler institute during the Spring 2005 semester on
``Algebraic combinatorics''.

\bibliographystyle{amsplain}

\providecommand{\bysame}{\leavevmode\hbox to3em{\hrulefill}\thinspace}
\providecommand{\MR}{\relax\ifhmode\unskip\space\fi MR }
\providecommand{\MRhref}[2]{%
  \href{http://www.ams.org/mathscinet-getitem?mr=#1}{#2}
}
\providecommand{\href}[2]{#2}

\end{document}